\documentclass[english,a4paper, 12pt,dvipdfm]{amsart}
\makeindex
\usepackage[all,web]{xy}
\usepackage{graphicx}
\usepackage{amsmath} 
\usepackage{amsfonts}
\usepackage{amssymb}
\usepackage{amsthm}
\usepackage{setspace}
\usepackage{geometry}
 \theoremstyle{plain}    
 \newtheorem{thm}{Theorem}[section]
 \numberwithin{equation}{section} 
 \numberwithin{figure}{section} 
 \theoremstyle{remark}   
  
 \theoremstyle{plain}    
 \theoremstyle{plain}    
  
 \theoremstyle{remark} 
 
 \theoremstyle{remark}
 
 \theoremstyle{definition}
  
 \theoremstyle{plain}
 \newtheorem{conjecture}[thm]{Conjecture} 
 \theoremstyle{plain}    
 \theoremstyle{plain}    
 \theoremstyle{definition}
 
 \theoremstyle{definition}
  \newtheorem{example}[thm]{Example}
 \theoremstyle{plain}    
  
 \theoremstyle{plain}    
 \newtheorem{lem}[thm]{Lemma} 
 \theoremstyle{remark}    
  
 \theoremstyle{remark}    
  
 \theoremstyle{definition}
  
 \theoremstyle{plain}
 \theoremstyle{remark}
 \newtheorem{rem}[thm]{Remark}
 \theoremstyle{remark}

\font\cyr=wncyr10 scaled \magstep 1
\newcommand{\Sha}{\mbox{\cyr X}}

\def\CC{\mathbb{C}}
\def\FF{\mathbb{F}}

\def\QQ{\mathbb{Q}}

\def\ZZ{\mathbb{Z}}

\def\inner<#1,#2>{{\left\langle{{#1},{#2}}\right\rangle}}
\def\sl2of#1{\textrm{SL}_2(#1)}

\def\inner<#1,#2>{{\left\langle{{#1},{#2}}\right\rangle}}

\def\ssp{\def\baselinestretch{1.0}\large\normalsize}

\usepackage{babel}

\begin{document}

\title{On the 2-divisibility of certain Heenger points}
\author{Carlos Casta\~no-Bernard}
\email{ccb1002@dpmms.cam.ac.uk}

\begin{abstract}
Let $E$ be an elliptic curve defined over the rationals and
let $N$ be its conductor. Assume $N$ is prime. In this paper we
give numerical evidence that suggests some conjectures on
the $2$-divisibility of certain sums of Heenger points 
of discriminant $D$ dividing $4N$ on the elliptic curve $E$.
One of these conjectures suggests a possible link between
the parity of the eigenvalue $a_A(2)$ and
the parity of the  \v Safarevi\v c-Tate group $\Sha(A)$ of
certain elliptic curves $A$ of square conductor.
\end{abstract}


\maketitle

\pagenumbering{roman}
\setcounter{page}{0}



\tableofcontents


\pagenumbering{arabic}
\pagestyle{headings}

\section{Preliminaries}
Let $E$ is an elliptic curve over $\QQ$.
Using Tate's algorithm~\cite{tate:algorithm}
we may compute a global
minimal Weierstra\ss\ equation
\begin{equation}\label{eqn:weier}
Y^2+a_1XY+a_3Y=X^3+a_2X^2+a_4X+a_6
\end{equation}
the conductor $N_E$ of $E$,
the local Tamagawa numbers $c_p$ for $p|N$,
and the minimal discriminant $\Delta_E$.
We assume $N_E$ odd to simplify the discussion.
The $L$-function of $E$ is given by
\begin{displaymath}
L(E,s)=\sum_{n=1}^\infty a_E(n)n^{-s}=
\prod_{p\,\textit{prime}}
(1-a_E(p)p^{-s} + \delta_p p^{1-2s})^{-1},
\end{displaymath}
where $a_E(p)=p+\delta_p-\#(E_{\FF_p})$,
with $\delta_p = 0$ if $E$ has bad
reduction at $p$ and $1$ otherwise.
By the work of Wiles~\cite{wiles:flt}
and
Breuil-Conrad-Diamond-Taylor~\cite{breuil:modularity}
we know that 
\begin{equation}
f(\tau)=\sum_{n=1}^\infty a_E(n)q^n\qquad\qquad (q=e^{2\pi \tau})
\end{equation}
is a normalized newform of weigth $2$ and level $N_E$.
Let $E_f$ be the elliptic curve over the field of
complex number $\CC$
determined by the periods of
\begin{equation}
\omega_f=2\pi i f(\tau)d\tau,
\end{equation}
regarding $\omega_f$ as a (holomorphic) differencial
on the modular curve $X_0(N)$.
By the work of
Mazur and Swinnerton-Dyer~\cite{mazur:weilcurves}
we know that $E_f$ may be provided with
a $\QQ$-structure such that
$L(E_f,s)=L(f,s)$
and such that the natural map
from $X_0(N)$ to $E_f$
sending $\infty$ to the identity element $O$ of $E_f$
is defined over $\QQ$.
In particular,
by Faltings' isogeny theorem~\cite{faltings:finiteness}
the curves $E$ and $E_f$ are isogenous over $\QQ$.
Thus there is a non-constant morphism
\begin{equation}
\varphi:X_0(N)\longrightarrow E.
\end{equation}
defined over $\QQ$.
Unless otherwise stated,
let us assume
the sign in the functional equation
$\Lambda(E,s)=\pm\Lambda(E,2-s)$ is ``$-$'',
where 
\begin{equation}
\Lambda(E,s)=N_E^{\frac{s}{2}}(2\pi)^{-s}\Gamma(s)L(E,s).
\end{equation}
In other words
$\varphi$ factors through the
Atkin-Lehner quotient $X_0^+(N)$
defined as $X_0(N)$ modulo the involution $w_N$.
Now suppose we have a pair $(D,r)$
that satisfies the so-called 
\textit{Heegner condition}\footnote{This condition was 
introduced by Birch~\cite{birch:heegner}.}
i.e. $D$ is the discriminant of
an imaginary quadratic order
$\mathcal{O}_D$
of conductor $f$ such that gcd$(N,f)=1$
and $r\in\ZZ$ is such that
\begin{equation}
D\equiv r^2\pmod{4N}.
\end{equation}
So we have a proper $\mathcal{O}_D$-ideal
$\mathfrak{N}=\ZZ N + \ZZ\frac{-r+\sqrt{D}}{2}\subset\mathcal{O}_D$
and $\mathcal{O}_D/\mathfrak{N}\cong \ZZ/N\ZZ$.
Using Gross's notation as in~\cite{gross:heegner},
for each proper $\mathcal{O}_D$-ideal $\mathfrak{A}$
consider the
\textit{Heegner point} $x=(\mathcal{O}_D,\mathfrak{N},[\mathfrak{A}])$,
which lies in $X_0(N)_H$,
where $H$ is the ring class field attached to $\mathcal{O}_D$,
and $[\mathfrak{A}]$ denotes the 
class of $\mathfrak{A}$ in Pic$(\mathcal{O}_D)$.
To simplify our exposition let us assume
$E_\QQ\cong \ZZ$.
Define the \textit{weighted trace} $t_{D,r,\varphi}$ of $\varphi(x)$ on $E$
by the equation
\begin{equation}
u_Dt_{D,r,\varphi}=
\sum_{\sigma\in\textit{Gal}(H/K)}\varphi (x^\sigma)\in E_\QQ,
\end{equation}
where
\begin{displaymath}
u_D=\left\{
\begin{array}{ll}
\frac{1}{2}\#(O^\times_D),&\textrm{if $\#(O^\times_D)>2$,}\\ 
2,&\textrm{if $\#(O^\times_D)=2$ and $N|D$,}\\
1,&\textrm{otherwise.}\\
\end{array}
\right.
\end{displaymath}
and $K=\QQ(\sqrt{D})$.
Now define the \textit{generalized trace} $y_{D,r,\varphi}\in E_\QQ$
as the sum
\begin{equation}
y_{D,r,\varphi}=\sum_{e^2|D}t_{\frac{D}{e^2},\frac{r}{e},\varphi}
\end{equation}
where we define $t_{\frac{D}{e^2},\frac{r}{e},\varphi}=O_E$
if $(\frac{D}{e^2},\frac{r}{e})$ does not satisfy the
Heegner condition.
Assume further that the derivative $L^\prime(E,1)\not=0$,
so that $E$ has analytic $1$.
By the work of Kolyvagin~\cite{kolyvagin:mordell}
we know that the rank of $E$ is $1$.
To simplify the exposition assume $E_\QQ$ is torsion free.
Fix a generator $g_E$ of the Mordell-Weil group $E_\QQ$ of $E$ over $\QQ$.
Define $\beta_{D,r}\in\ZZ$ by the equation
\begin{equation*}
(y_{D,r}^+)_f=\beta_{D,r}g_E,
\end{equation*}
where $(y_{D,r}^+)_f$ denotes the $f$-eigencomponent of $y_{D,r}^+$,
and $y_{D,r}^+$ denotes the image of $y_{D,r}$ in the Jacobian $J_0^+(N)$
as in Gross, Kohnen and Zagier~\cite{gross:gkz}.
Since $w_N$ acts as complex conjugation on $y_{D,r}$,
we see $y_{D,r}^+$ is defined over $\QQ$,
and the definition of $\beta_{D,r}$ makes sense.
To ease notation we write $\beta_D$ instead of $\beta_{D,r}$
if there is no risk of confusion.

%
\section{The conjectures}
Suppose $N$ divides $D$.
So taking the weighted trace associated to $D$ involves
dividing by $2$ a sum over a Galois orbit over $\QQ$ of points on $E$
obtained as images of the Heegner points of discriminant $D$.
Now consider in particular
the discriminant $D=-4N$ and assume from now on that $N$ is prime.
Numerical evidence shown below strongly suggests that
(at least) the $f$-eigencomponents $(y_D^+)_f$
of the generalised trace $y_D$
may be further divided by $2$ in $E(\QQ)$,
provided $N$ satisfies a certain
congruence condition.\footnote{A consecuence of
Proposition 3.1 (p. 347) of Gross~\cite{gross:primelevel} is
that the Heenger points involved in the generalised trace
of discriminant $D=-4N$ are precisely
the fixed points of the Fricke involution $w_N$.
But we do not use this fact here.}
But we shall find convenient to
state a slightly stronger conjecture first.

\begin{table}
\ssp
\begin{center}
\caption{$\beta_{-N,E}$ such that $a_E(2)\in 2\ZZ$
and $N\equiv 3\pmod{4}$ prime.\label{tbl:a2even}}
\end{center}%
\vspace{-.3in}%
$$
\begin{array}{lrr}
E & a_E(2)& \beta_{-N}\\
\vspace{-2ex}\\
\mathbf{43A} &-2 & 2\\
\mathbf{ 131A} &0 & 0\\
\mathbf{ 163A} &0 & -2\\
\mathbf{ 347A} &-2 & -2\\
\mathbf{ 443A} &0 & -2\\	
\vspace{-2ex}\\
\mathbf{ 467A} &0 & 0\\
\mathbf{ 811A} &0 & 2\\
\mathbf{ 827A} &0 & 2\\
\mathbf{ 1019A} &0 & 0\\
\mathbf{ 1019B} &-2 & -4\\	
\vspace{-2ex}\\
\mathbf{ 1051A} &0 & 0\\
\mathbf{ 1259A} &0 & -2\\
\mathbf{ 1747A} &2 & 2\\
\mathbf{ 1811A} &0 & -2\\
\mathbf{ 1987A} &0 & -2\\	
\vspace{-2ex}\\
\mathbf{ 2539A} &-2 & 4\\
\mathbf{ 2699A} &0 & 0\\
\mathbf{ 3251A} &0 & 2\\
\mathbf{ 3259A} &-2 & -4\\
\mathbf{ 3259B} &-2 & -10\\	
\vspace{-2ex}\\
\mathbf{ 3347A} &2 & 0\\
\mathbf{ 3547A} &0 & -2\\
\mathbf{ 3851A} &-2 & 0\\
\mathbf{ 3931A} &0 & -2\\
\mathbf{ 3947A} &0 & 0\\	
\vspace{-2ex}\\
\mathbf{ 4051A} &0 & 4\\
\mathbf{ 4507A} &-2 & -6\\
\mathbf{ 4603A} &-2 & -2\\
\mathbf{ 5443A} &2 & 0\\
\mathbf{ 5563A} &0 & 4\\	
\vspace{-2ex}\\
\mathbf{ 6131A} &2 & 0\\
\mathbf{ 6691A} &0 & -8\\
\mathbf{ 7019A} &-2 & -2\\
\mathbf{ 7187A} &2 & 2\\
\mathbf{ 7283A} &0 & 0\\	
\end{array}\qquad\qquad
\begin{array}{lrr}
 E & a_E(2) & \beta_{-N}\\
\vspace{-2ex}\\
\mathbf{ 8419A} &2 & 0\\
\mathbf{ 8747A} &0 & -4\\
\mathbf{ 8803A} &0 & -2\\
\mathbf{ 9539A} &2 & 0\\
\mathbf{ 9587A} &2 & -2\\	
\vspace{-2ex}\\
\mathbf{ 9811A} &0 & -4\\
\mathbf{ 10859A} &0 & 4\\
\mathbf{ 10859B} &-2 & 4\\
\mathbf{ 10987A} &0 & 0\\
\mathbf{ 11867A} &-2 & 6\\	
\vspace{-2ex}\\
\mathbf{ 11923A} &0 & -4\\
\mathbf{ 11939A} &0 & -2\\
\mathbf{ 11939B} &2 & 0\\
\mathbf{ 12163A} &2 & 4\\
\mathbf{ 12619A} &-2 & -6\\	
\vspace{-2ex}\\
\mathbf{ 13043A} &0 & 0\\
\mathbf{ 13523A} &2 & 0\\
\mathbf{ 15083A} &-2 & 6\\
\mathbf{ 15091A} &2 & 2\\
\mathbf{ 15131A} &0 & 4\\	
\vspace{-2ex}\\
\mathbf{ 15227A} &0 & -2\\
\mathbf{ 15971A} &2 & 2\\
\mathbf{ 16883A} &0 & 2\\
\mathbf{ 16963A} &2 & 10\\
\mathbf{ 17387A} &0 & 0\\	
\vspace{-2ex}\\
\mathbf{ 17387B} &-2 & -8\\
\mathbf{ 17483A} &-2 & 6\\
\mathbf{ 17747A} &2 & 2\\
\mathbf{ 17827A} &0 & 4\\
\mathbf{ 18059A} &0 & 0\\	
\vspace{-2ex}\\
\mathbf{ 18251A} &-2 & 8\\
\mathbf{ 18859A} &-2 & -2\\
\mathbf{ 19387A} &0 & -4\\
\mathbf{ 19387B} &0 & -6\\
&&\\
\end{array}
$$
\end{table}

\begin{conjecture}\label{conj:a2even}
Assume that $E$ is such that $N_E\equiv 3\pmod{4}$.
If $a_E(2)$ is even then $\beta_{-N}$ is even.
(See Table~\ref{tbl:a2even}, below.)
\end{conjecture}

By the work of Gross-Kohnen-Zagier~\cite{gross:gkz} 
we know that for each prime $p$ such that gcd$(p,N)=1$ we have
\begin{equation}\label{eqn:jacobihecke}
\beta_{p^2D,pr,\varphi}+\left(\frac{D}{p}\right)\beta_{D,r,\varphi}=
a_E(p)\beta_{D,r,\varphi},
\end{equation}
provided $D$ is fundamental.
So clearly Equation~\ref{eqn:jacobihecke},
with $p=2$ and $D=-N$ together with
Conjecture~\ref{conj:a2even} imply the following.

\begin{conjecture}\label{conj:minus4n}
Assume that $E$ is such that $N_E\equiv 3\pmod{4}$.
Then the integer $\beta_{-4N}$ is even.
\end{conjecture}

\begin{example}\label{ex:43a1}
Suppose $E$ is elliptic curve \textbf{43A1}
of Cremona's Tables~\cite{cremona:onlinetables}.
So $E$ is an elliptic curve with minimal Weierstra\ss\ equation
\begin{displaymath}
Y^2 + Y = X^3 + X^2,
\end{displaymath}
$E_\QQ=\ZZ g_E$, where $g_E=(0,0)$, 
and the eigenvalue $a_E(2)=-2$.
Note that the pair $(-43,129)$
satisfies the Heegner condition.
We may find $\beta_{-43,129}=(-1,-1)\in E(\QQ)$.
So $y_{-N}=2g_E$,
and thus $\beta_{-N}=2$, which is even.
So Conjecture~\ref{conj:a2even} holds for curve \textbf{43A1}.
Now using Equation~\ref{eqn:jacobihecke} we
have $\beta_{-4\cdot43}=\left(a_E(2)-\left(\frac{-43}{2}\right)\right)\beta_{-43}$.
So $\beta_{-4\dot 43}=(-2+1)\beta_{-43}=-\beta_{-43}=-2g_E$.
In other words $\beta_{-4\cdot 43}=-2$ and
Conjecture~\ref{conj:minus4n} holds for curve \textbf{43A1}.
\end{example}

Using the following lemma we will actually prove
part of Conjecture~\ref{conj:minus4n}.

\begin{lem}\label{lem:7mod8}
Let $E$ be an elliptic curve of prime conductor $N$.
If $N\equiv 7\pmod{8}$,
then $a(2)=\pm 1$. 
\end{lem}

\begin{proof}
By brute force it may be found the list of the reduction modulo $p=2$
of all possible Weierstra\ss\ models of $E$ over $\ZZ$
together with the corresponding eigenvalue $a(p)$:
\begin{displaymath}
\begin{array}{rrrrrr}
a_1 & a_2 & a_3 & a_4 & a_6& \, a(2)\\
\vspace{-2ex}\\
0& 0& 1& 0& 0 & \, 0\\
0& 0& 1& 0& 1 & \, 0\\
0& 0& 1& 1& 0 & \, -2\\
0& 0& 1& 1& 1 & \, 2\\
\vspace{-2ex}\\
0& 1& 1& 0& 0 & \, -2\\
0& 1& 1& 0& 1 & \, 2\\
0& 1& 1& 1& 0 & \, 0\\
0& 1& 1& 1& 1 & \, 0\\
\end{array}\qquad\qquad
\begin{array}{rrrrrr}
a_1 & a_2 & a_3 & a_4 & a_6& \, a(2)\\
\vspace{-2ex}\\
1& 0& 0& 0& 1 & \, -1\\
1& 0& 0& 1& 0 & \, -1\\
1& 0& 1& 0& 1 & \, 1\\
1& 0& 1& 1& 1 & \, 1\\
\vspace{-2ex}\\
1& 1& 0& 0& 1 & \, 1\\
1& 1& 0& 1& 0 & \, 1\\
1& 1& 1& 0& 0 & \, -1\\
1& 1& 1& 1& 0 & \, -1\\
\end{array}
\end{displaymath}
By inspection we may see that $a(2)\in 2\ZZ$ implies $a_1\in 2\ZZ$.
Now suppose $(a_1,a_2,a_3,a_4,a_6)$
is a global minimal Weierstra\ss\ equation of $E$ over $\QQ$,
normalised so that $|a_3|\leq 1$.
A straight forward computation shows that
\begin{displaymath}
\begin{split}
\Delta \equiv 7a_6a_1^6 + a_4a_3a_1^5 + 7a_3^2a_2a_1^4 + 4a_6a_2a_1^4 +
a_4^2a_1^4 +\\
a_3^3a_1^3 + 4a_6a_3a_1^3 + 2a_4a_3^2a_1^2 + 4a_3^3a_2a_1 + 5a_3^4\pmod{8},
\end{split}
\end{displaymath}
where $\Delta$ is the discriminant of the Weierstra\ss\ equation.
So $a_1\in 2\ZZ$ implies $\Delta\equiv 5a_3^4\pmod{8}$.
But we assumed $N$ prime.
Therefore $\Delta=\pm N$ and $a_3=\pm 1$.
Thus $N\equiv\pm 5\pmod{8}$,
i.e.  $N\not\equiv 7\pmod{8}$ and $a(2)$ is odd.
In other words $a(2)=\pm 1$
and the proof of the lemma is now complete.
\end{proof}

\begin{thm}
If $E$ is such that $N\equiv 7\pmod{8}$,
then $\beta_{-4N}$ is even.
\end{thm}

\begin{proof}
Since $N$ is prime and $N\equiv 7\pmod{8}$,
Lemma~\ref{lem:7mod8} implies $a_E(2)$ is odd.
So if we set $p=2$ and $D=-N$ in Equation~\ref{eqn:jacobihecke}
we see $\beta_{-4N,E}$ is even.
\end{proof}

\begin{rem}
The entries of
Table~\ref{tbl:a2odd} are again the numbers $\beta_{-N}$,
for prime $N\leq 20\,000$ such that $N\equiv 3\pmod{4}$,
now with  $a_E(2)$ \textit{odd}.
In sharp contrast with Table~\ref{tbl:a2even},
in Table~\ref{tbl:a2odd} we may see examples of both,
even $\beta_{-N}$ and odd $\beta_{-N}$.
There are $25$ even $\beta_{-N}$ and $33$ odd $\beta_{-N}$.
\end{rem}

\begin{table}
\ssp
\begin{center}
\caption{$\beta_{-N,E}$ such that $a_E(2)=\pm 1$ and
$N\equiv 3\pmod{4}$ prime.\label{tbl:a2odd}}
\end{center}%
\vspace{-.3in}%
$$
\begin{array}{lrr}
E & a_E(2)& \beta_{-N}\\
\vspace{-2ex}\\
\mathbf{ 79A} &-1 & 1\\
\mathbf{ 83A} &-1 & 1\\
\mathbf{ 331A} &-1 & 4\\
\mathbf{ 359A} &1 & 2\\
\mathbf{ 359B} &-1 & 0\\	
\vspace{-2ex}\\
\mathbf{ 431A} &-1 & 1\\
\mathbf{ 443B} &-1 & 1\\
\mathbf{ 503A} &1 & 1\\
\mathbf{ 659A} &1 & 1\\
\mathbf{ 1091A} &-1 & -3\\	
\vspace{-2ex}\\
\mathbf{ 1439A} &1 & 1\\
\mathbf{ 1607A} &-1 & -1\\
\mathbf{ 3023A} &-1 & 1\\
\mathbf{ 3163A} &1 & -4\\
\mathbf{ 3391A} &-1 & 1\\	
\vspace{-2ex}\\
\mathbf{ 3803A} &1 & 3\\
\mathbf{ 4159A} &1 & 0\\
\mathbf{ 4159B} &1 & -6\\
\mathbf{ 4799A} &-1 & 2\\
\mathbf{ 4799B} &-1 & 6\\	
\vspace{-2ex}\\
\mathbf{ 5503A} &-1 & -3\\
\mathbf{ 5867A} &1 & 0\\
\mathbf{ 5987A} &1 & 2\\
\mathbf{ 6011A} &1 & 1\\
\mathbf{ 6199A} &-1 & -4\\	
\vspace{-2ex}\\
\mathbf{ 6427A} &1 & 1\\
\mathbf{ 6823A} &-1 & -1\\
\mathbf{ 6967A} &1 & 0\\
\mathbf{ 7219A} &-1 & -1\\
\mathbf{ 7699A} &1 & 2\\	
\end{array}\qquad\qquad
\begin{array}{lrr}
E & a_E(2) & \beta_{-N}\\
\vspace{-2ex}\\
\mathbf{ 7723A} &1 & 5\\
\mathbf{ 8167A} &1 & -1\\
\mathbf{ 8623A} &1 & 3\\
\mathbf{ 9127A} &-1 & -4\\
\mathbf{ 9491A} &-1 & 4\\	
\vspace{-2ex}\\
\mathbf{ 9811B} &-1 & -1\\
\mathbf{ 10163A} &1 & 0\\
\mathbf{ 10567A} &1 & -3\\
\mathbf{ 10799A} &1 & 5\\
\mathbf{ 11119A} &1 & -1\\	
\vspace{-2ex}\\
\mathbf{ 12007A} &1 & 13\\
\mathbf{ 12227A} &1 & -5\\
\mathbf{ 12547A} &1 & -2\\
\mathbf{ 13451A} &-1 & -7\\
\mathbf{ 13619A} &-1 & -2\\	
\vspace{-2ex}\\
\mathbf{ 13723A} &-1 & 10\\
\mathbf{ 13723B} &-1 & 4\\
\mathbf{ 15551A} &1 & 1\\
\mathbf{ 15859A} &-1 & 1\\
\mathbf{ 16411A} &1 & 4\\	
\vspace{-2ex}\\
\mathbf{ 17299A} &-1 & -6\\
\mathbf{ 18059B} &1 & 2\\
\mathbf{ 18127A} &1 & -8\\
\mathbf{ 18523A} &-1 & -7\\
\mathbf{ 18899A} &1 & -5\\	
\vspace{-2ex}\\
\mathbf{ 19211A} &1 & -6\\
\mathbf{ 19583A} &-1 & -1\\
\mathbf{ 19927A} &1 & -3\\
&&\\
&&\\
\end{array}$$
\end{table}

Let $\Sha(E)$ be the \v Safarevi\v c-Tate group of $E$ over $\QQ$
and let $E^D$ be the twist of the elliptic curve $E$
associated to a fundamental discriminant $D$.
Conjecture~\ref{conj:a2even} suggests the following.

\begin{conjecture}\label{conj:squarelevel}
Let $E$ be as above and assume $N\equiv 3\pmod{4}$,
with $N$ prime such that $\beta_{-N}\not=0$.
If $a_E(2)$ is even then $|\Sha(E^{-N})|$ is even.
\end{conjecture}

Recall that for each fundamental discriminant $D$ and each prime $p$
\begin{equation}\label{eqn:twisteigen}
a_A(p)=\chi(p)a_E(p),
\end{equation}
where $\chi$ is the quadratic character attached to $\QQ(\sqrt{D})$
and $A=E^D$.
In particular,
by putting $p=2$ and $D=-N$ in
Equation~\ref{eqn:twisteigen} (and noting $\chi(2)\not=0$)
we may see $a_A(2)$ and $a_E(2)$ have the same parity.
Also note $N_A=N_E^2$, where $N_A$ denotes the conductor of $A$. 
So Conjecture~\ref{conj:squarelevel} suggests a possible
link between
the parity of $a_A(2)$ and the parity of $|\Sha(A)|$,
for certain elliptic curves $A$ such that $N_A$ is square.

\begin{example}
Suppose $E$ is curve {\bf 43A1}.
(See Example~\ref{ex:43a1},
above.)
With the help of Tate's algorithm~\cite{tate:algorithm}
it is easy to find a global minimal model for $A=E^{(-N)}$.
Then using Cremona's Tables~\cite{cremona:onlinetables}
it may be found that $A$ is in fact
curve {\bf 1849D1} of conductor $N_A=43^2$.
A short calculation shows that the eigenvalue $a_A(2)=2$.
This is consistent with the (refined version of the)
Birch and Swinnerton-Dyer conjecture,
which predicts that $|\Sha(A)|=4$.
\end{example}

\begin{table}
\ssp
\begin{center}
\caption{$\beta_{-4N,E}$ and $\beta_{-4}$ with $N\equiv 1\pmod{4}$ prime.
\label{tbl:1mod4}}
\end{center}%
\vspace{-.3in}%
$$
\begin{array}{lrr}
E & \beta_{-4}& \beta_{-4N}\\
\vspace{-2ex}\\
{\bf 37A} &-1 & 3\\
{\bf 53A} &-1 & 1\\
{\bf 61A} &1 & 1\\
{\bf 89A} &1 & -1\\
{\bf 101A} &-1 & 1\\	
\vspace{-2ex}\\
{\bf 197A} &-1 & -5\\
{\bf 229A} &0 & -2\\
{\bf 269A} &-1 & 3\\
{\bf 277A} &-1 & 3\\
{\bf 373A} &1 & 1\\	
\vspace{-2ex}\\
{\bf 557A} &1 & 1\\
{\bf 593A} &1 & -5\\
{\bf 677A} &-1 & -1\\
{\bf 797A} &0 & 0\\
{\bf 829A} &-1 & -9\\	
\vspace{-2ex}\\
{\bf 997A} &0 & 8\\
{\bf 1549A} &1 & 5\\
{\bf 1949A} &1 & 3\\
{\bf 1973A} &1 & 1\\
{\bf 2017A} &1 & 5\\	
\vspace{-2ex}\\
{\bf 2089A} &0 & 8\\
{\bf 2141A} &-1 & -3\\
{\bf 2161A} &1 & 1\\
{\bf 2221A} &1 & 1\\
{\bf 2269A} &1 & -3\\	
\vspace{-2ex}\\
{\bf 2341A} &-1 & 1\\
{\bf 2357A} &1 & -1\\
{\bf 2557A} &0 & 2\\
{\bf 2609A} &1 & 1\\
{\bf 2749A} &1 & -5\\	
\vspace{-2ex}\\
{\bf 3109A} &1 & 7\\
{\bf 3229A} &0 & 6\\
{\bf 3313A} &-1 & 3\\
{\bf 3469A} &-1 & -7\\
{\bf 3797A} &-2 & -10\\	
\end{array}\quad
\begin{array}{lrr}
E & \beta_{-4} & b_{-4N}\\
\vspace{-2ex}\\
{\bf 3853A} &1 & 3\\
{\bf 3877A} &1 & 5\\
{\bf 4021A} &1 & -1\\
{\bf 4481A} &-2 & -10\\
{\bf 4481B} &0 & 4\\	
\vspace{-2ex}\\
{\bf 4493A} &-2 & 4\\
{\bf 5237A} &-1 & 5\\
{\bf 5309A} &-1 & 7\\
{\bf 5417A} &1 & 3\\
{\bf 5417B} &4 & -8\\	
\vspace{-2ex}\\
{\bf 5417C} &0 & -4\\
{\bf 5653A} &0 & 2\\
{\bf 5717A} &0 & -2\\
{\bf 6373A} &-2 & 10\\
{\bf 6689A} &1 & -5\\	
\vspace{-2ex}\\
{\bf 7109A} &1 & -1\\
{\bf 7109B} &1 & 3\\
{\bf 7213A} &1 & 15\\
{\bf 7757A} &1 & 3\\
{\bf 8069A} &0 & -4\\	
\vspace{-2ex}\\
{\bf 8101A} &-3 & -19\\
{\bf 8597A} &0 & -4\\
{\bf 8929A} &1 & 7\\
{\bf 9109A} &1 & 3\\
{\bf 9413A} &0 & -10\\	
\vspace{-2ex}\\
{\bf 9829A} &-4 & -4\\
{\bf 9941A} &-1 & -1\\
{\bf 10061A} &-1 & -3\\
{\bf 10333A} &1 & -5\\
{\bf 10333B} &1 & 7\\	
\vspace{-2ex}\\
{\bf 10733A} &0 & -2\\
{\bf 10789A} &-1 & -1\\
{\bf 10949A} &2 & -4\\
{\bf 11321A} &1 & -1\\
{\bf 11321B} &0 & -4\\	
\end{array}\quad
\begin{array}{lrr}
E & \beta_{-4} & b_{-4N}\\
\vspace{-2ex}\\
{\bf 11353A} &1 & -5\\
{\bf 11789A} &0 & 10\\
{\bf 12097A} &-1 & 1\\
{\bf 12277A} &-1 & -1\\
{\bf 12289A} &-1 & -7\\	
\vspace{-2ex}\\
{\bf 12413A} &1 & 3\\
{\bf 13093A} &0 & -6\\
{\bf 13537A} &2 & 14\\
{\bf 13789A} &2 & 14\\
{\bf 14173A} &-1 & 1\\	
\vspace{-2ex}\\
{\bf 14461A} &0 & -8\\
{\bf 15013A} &-1 & 9\\
{\bf 15101A} &-1 & 15\\
{\bf 15349A} &0 & -2\\
{\bf 15641A} &1 & -5\\	
\vspace{-2ex}\\
{\bf 15661A} &1 & -31\\
{\bf 15661B} &1 & -5\\
{\bf 15773A} &0 & -6\\
{\bf 15889A} &-1 & 5\\
{\bf 16061A} &-1 & 1\\	
\vspace{-2ex}\\
{\bf 16189A} &-1 & -11\\
{\bf 16369A} &0 & -8\\
{\bf 16649A} &0 & 4\\
{\bf 16649B} &0 & -4\\
{\bf 16889A} &-1 & 5\\	
\vspace{-2ex}\\
{\bf 16937A} &1 & 3\\
{\bf 17093A} &1 & -5\\
{\bf 17573A} &0 & 12\\
{\bf 17837A} &-1 & 9\\
{\bf 18097A} &0 & 12\\	
\vspace{-2ex}\\
{\bf 18269A} &0 & -2\\
{\bf 18397A} &0 & 12\\
{\bf 19469A} &2 & -4\\
&&\\
&&\\
\end{array}$$
\end{table}

We have a further conjecture.
Table~\ref{tbl:1mod4} suggests the following.

\begin{conjecture}\label{conj:1mod4}
If $N\equiv 1\pmod{4}$ then
the integers $\beta_{-4N}$ and $\beta_{-4}$
have the same parity.
(See Table~\ref{tbl:1mod4}, below.)
\end{conjecture}

%
\section{Further remarks}
A well-known result of Gross, Kohnen, and Zagier~\cite{gross:gkz}
implies that the integers $\beta_D$ for $D=-3,-4,\dots$ are
the coefficients of the newform $g$ of weight $\frac{3}{2}$,
that corresponds (modulo multiplication by a scalar)
to the newform $f$ of weight $2$ via
the Shimura correspondence.\footnote{Recall we assumed $N$ prime;
for $N$ composite we would need to consider instead
of newforms $g$ of weight $\frac{3}{2}$,
Jacobi newforms of weight $2$ and index $N$.}
Perhaps a divisibility theory of (suitably normalised)
half-integral weight modular forms
analogous to the divisibility theory of integral weight
modular forms might lead to a proof of the conjectures stated above. 
However,
not much is known about divisibility properties of
modular forms of half-integral weight,
apart from the work of Koblitz~\cite{koblitz:padiccon},
and
the work of Balog, Darmon, and Ono~\cite{darmon:congruences},
which is about modular forms of half-integral weight and level $N=4$.
For example,
by the work of Koblitz~\cite{koblitz:padiccon}
we know that Ramanujan's famous congruence
\begin{equation*}
\tau(n)\equiv \sigma_{11}(n)\pmod{691},
\end{equation*}
where 
\begin{equation*}
\Delta=q\prod_{n=1}^\infty(1-q^n)^{24}=\sum_{n=1}^\infty\tau(n)q^n
\end{equation*}
and $\sigma_k(n)=\sum_{d|n}d^k$
descends via the Shimura correspondence to the
congruence
\begin{equation*}
\alpha(n)\equiv -252 c(n)\pmod{691},
\end{equation*}
where $\delta=\sum_{n=1}^\infty \alpha(n)q^n$
is the cusp form of weight $\frac{13}{2}$ for $\Gamma_0(4)$
that corresponds to $\Delta$ under the Shimura lift,
normalised so that $\alpha(1)=1$
and $c(n)$ is the coefficient of $q^n$ in H. Cohen's
generalised class number
Eisenstein series $H_{\frac{13}{2}}=\sum_{n=1}^\infty c(n)q^n$.

\bibliographystyle{amsplain}
\bibliography{biblio}
\printindex
\end{document}